\documentclass[12pt]{article}
\usepackage{graphicx}
\usepackage{amsfonts}

\newcommand{\pic}[4]{\vspace{1ex}\setlength{\unitlength}{1cm}
\begin{picture}(0,#3)(5,.5)
\put(#2){\includegraphics[#4]{#1.eps}}
\end{picture}\vspace{1ex}}

\newtheorem{proposition}{Proposition}[section]
\newtheorem{coro}[proposition]{Corollary}

\newtheorem{lemma}[proposition]{Lemma}
\def\proof{\noindent\textit{Proof. }}

\def\C{\mathbb{C}}
\def\cg{\overline} 

\def\D{\mathbb{D}}
\def\DS{\displaystyle}

\def\R{\mathbb{R}}
\def\S{{\mathcal{S}}}
\def\qed{\hskip1em\raise3.5pt\hbox{\framebox[2mm]{\ }}}
\def\re{\,\!\mbox{Re}\,}
\def\im{\,\!\mbox{Im}\,}
\def\scq{s.c.q.}
\def\sgn{\,\!{\rm sgn}\,}

\newcommand{\X}[1]{X^{(#1)}}
\newcommand{\Xt}[1]{\widetilde X^{(#1)}}

\begin{document}
  
\begin{center} 
{\LARGE 
 Conformal Mapping  of Circular Quadrilaterals and
Weierstrass Elliptic Functions \\[2ex] \Large
  Philip R.  Brown\\
R. Michael Porter\footnote{Partially supported by CONACyT 80503}
}
  
\bigskip
\end{center}

\begin{minipage}{\textwidth}

\small 

\noindent Abstract.  Numerical and theoretical aspects of conformal
mappings from a disk to a circular-arc quadrilateral, symmetric with
respect to the coordinate axes, are developed.  The problem of
relating the accessory parameters (prevertices together with
coefficients in the Schwarzian derivative) to the geometric parameters
is solved numerically, including the determination of the parameters
for univalence. The study involves the related mapping from an
appropriate Euclidean rectangle to the circular-arc quadrilateral.
Its Schwarzian derivative involves the Weierstrass $\wp$-function, and
consideration of this related mapping problem leads to some new
formulas concerning the zeroes and the images of the half-periods of
$\wp$.



\medskip
   
\end{minipage}

\setcounter{section}{-1}
\section{Introduction \label{sec:intro}}

Conformal mappings of the unit disk onto a symmetric circular arc
quadrilateral (\scq) with right angles at the vertices were
investigated in \cite{Br1,Br2,KP2}.  The space of conformal
classes of such quadrilaterals is two-dimensional.  One parameter may
be taken to be the conformal module \cite{Ahl} of the quadrilateral,
and the other as a real factor $\lambda$ appearing in the Schwarzian
derivative of the conformal mapping, which is a rational function
whose general form is classical.  The articles cited are principally
concerned with the question of calculating these parameters numerically
when given the specific geometry of the quadrilateral, and the range
of values for which the Schwarzian derivative defines a univalent
function.

The aim of this article is to introduce an additional parameter in the
problem by allowing a variable angle of $\alpha\pi$ radians at each
vertex, where $0\le\alpha\le2$.  The numerical computations and the
univalence criteria are generalized to this context and refined.

An innovation here is to take into account the conformal mapping from
a rectangle (rather than the unit disk) to the \scq {} The advantage of
this approach is that the Schwarzian derivative of the mapping takes
the simple form $(1/2-2\alpha^2)\wp-2\mu$, and the Weierstrass
$\wp$-function has well known properties. The parameter $\mu \in \R$
plays the same role as the parameter $\lambda$ in \cite{KP2}: in rough
terms, either of these parameters controls the bending (i.e.,
curvatures) of the circular arcs bounding the \scq {} We work out
the explicit relationship between $\lambda$ and $\mu$ by making use of
the transformation (elliptic integral) of the unit disk onto a
rectangle. Let $\pm \exp{(\pm it)}$ be the prevertices in the unit
circle which map to the vertices of the rectangle; it is perhaps
surprising that the expression for $\mu$ in terms of $\lambda$ depends
on both $t$ and $\alpha$ (see (\ref{eq:mu}) below).

Our use of the $\wp$-function to clarify some aspects of the conformal
mapping problem of \scq s happily also sheds some light on the nature
of the $\wp$-function itself.  It is interesting, and apparently not
well known, that the constants $e_1$, $e_2$, and $e_3$ related to the
$\wp$-function have simple expressions in terms of trigonometric
ratios of $t$ (see (\ref{eq:e1e2e3}) below). Further, we make some
observations (in \ref{subsec:zeros}) on the location of the zeros of
$\wp$.
 
It is well known that the conformal mapping of the unit disk onto a
specified circular arc polygon can be expressed as a ratio of linearly
independent solutions of a second order ordinary differential
equation, with some parameters that have to be determined if the
polygon has more than three boundary arcs. If the prevertices and the
value of $\alpha$ are fixed, then the conformal mapping onto an
\scq\ has one unknown parameter. In \cite{KP2} the second order
differential equation is expressed as a Sturm-Liouville differential
equation with parameter $\lambda$, and the spectral parameter power
series (SPPS) method is introduced to solve the differential
equation. Here we use the SPPS method with parameter $\mu$ for
computing the conformal mapping from a rectangle. A ratio of linearly
independent solutions of the differential equation is in general
(i.e., for any value of $\mu$) conformal locally, but need not be
univalent. In Section~\ref{sec:domuniv} below we show how the Sturm
comparison theorem can be used to find bounds for $\mu$ which are
necessary for univalence of the conformal mapping. This makes it
possible to find a good ``seed value'' for the SPPS method, and as
explained and demonstrated in \cite{KP2} it is then possible to find
numerically the exact bounds for $\mu$ for which the conformal mapping
is univalent.

As an application of the SPPS method (with $\alpha=0$) we compute the
Fuchsian group for the universal covering map of a twice punctured
disk. More precisly, for a unit disk punctured at points $-a$ and $a$
for any given value \mbox{$0<a<1$} we are able to find the M\"obius
transformations of the Fuchsian group by calculating the value of a
parameter $s$ related to $a$. This is described in detail in
Section~\ref{sec:dpd} below. As an application of the explicit 
calculation of conformal mappings, some plots are given of
geodesics surrounding the two punctures for different values of $a$.

\section{Description of problem \label{sec:quads}}

\begin{figure}[b!]
\pic{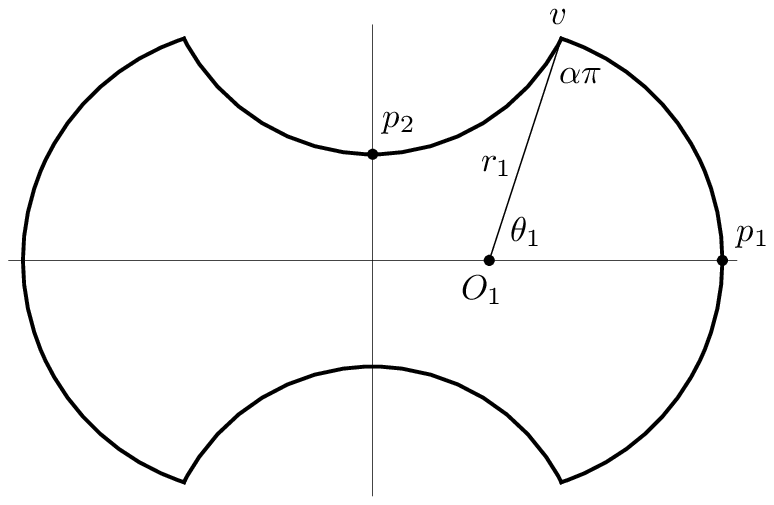}{8,.3}{4.5}{}

\pic{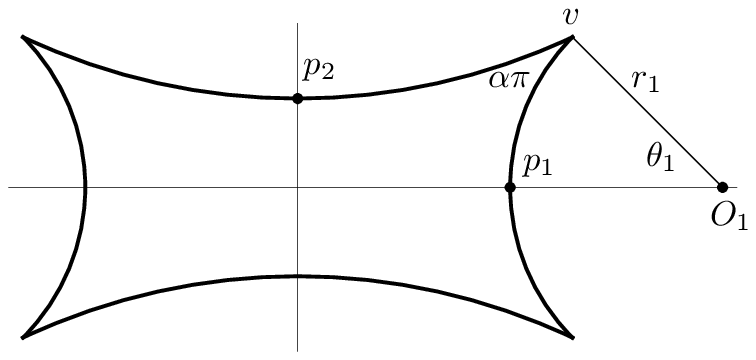}{8.77,0}{3.}{}
\caption{\small Geometric parameters defining symmetric circular 
 quadrilaterals.\label{fig:scq}}
\end{figure}

In this article a \textit{symmetric circular quadrilateral} (\scq) is
a domain $D$ in the complex plane $\C$ bounded by a Jordan curve
$\partial D$ which is formed joining four vertices $v$,
$\cg{v}$, $-v$, $-\cg{v}$ by circular arcs with common
internal angle $\alpha\pi$ at each vertex (Figure \ref{fig:scq}).  By
definition an \scq\ is symmetric in the real and imaginary axes.

\subsection{Parameters of \scq s \label{subsec:descrip}}

First we need to describe the geometry of an \scq\ $D$.  Denote by
$p_1$, $p_2$ the midpoints of the rightmost and the uppermost edges of
$D$.  Let $2\theta_1$, $2\theta_2$ be the angles subtended by these
edges with respect to their centers $O_1$, $O_2$, and let $r_1$, $r_2$
be the corresponding radii.  Necessarily
\begin{equation} 
   0 \le \alpha\pi \le 2\pi,\quad   0<\theta_i\le \pi.
\end{equation}

Mappings for $\alpha=1/2$ were studied in \cite{Br1,Br2,KP2}, and for
$\alpha=0$ in \cite{Phoro}.

Let $p_1$, $\theta_1$, and $r_1$ be given.  To find $v$ we need to
know whether the right edge is convex.  Let $\kappa_1$, $\kappa_2$
denote the signed curvatures of the edges, that is,
\[   \kappa_1 = \frac{1}{p_1-O_1}, \quad
     \kappa_2 = \frac{i}{p_2-O_2},
\]
so $r_i=1/|\kappa_i|$.  If $O_i=\infty$ then $\kappa_i=0$ while
$\theta_i$ is not defined.
 
Let $\alpha$ also be given.  Suppose first that $\kappa_1>0$. The
tangent direction of the right edge at $v$ leading towards $p_1$ is
$e^{i(\theta_1-\pi/2)}$ (necessarily pointing into the right
half-plane), so the direction of the upper edge at $v$ leading towards
$p_2$ is $e^{i(\theta_1-\pi/2-\alpha\pi)}$. It follows that the nature
of bending of the upper edge is given by
$\sgn\kappa_2=\sgn(\sin(\theta_1-\pi/2-\alpha\pi))$ (although we
cannot yet calculate $\kappa_2$).  If $\sgn\kappa_2<0$, then $O_2$
lies on the ray from $v$ in the direction of
$e^{i(\theta_1-\pi-\alpha\pi)}$, whereas if $\sgn\kappa_2>0$, the
direction is $e^{i(\theta_1-\alpha\pi)}$.  In either event this ray
meets the imaginary axis in $O_2=i(\im v -\re
v\tan(\theta_1-\alpha\pi))$, and finally the values $r_2=|v-O_2|$ and
then $\kappa_2=(\sgn\kappa_2)/r_2$ are determined.

Suppose on the other hand that $\kappa_1<0$. Then the tangent
direction of the right edge at $v$ leading towards $p_1$ is
$e^{-i(\theta_1+\pi/2)}$ (necessarily pointing into the left
half-plane), so the direction of the upper edge at $v$ leading towards
$p_2$ is now $e^{-i(\theta_1+\pi/2+\alpha\pi)}$.  Then $\sgn\kappa_2$
is the sign of the imaginary part of this number, and $O_2$ is
calculated similarly.  In summary, starting from   
\begin{eqnarray*}  
   r_1  &=& \frac{1}{|\kappa_1|},   \\
   O_1  &=&  p_1 - \frac{1}{\kappa_1}, \\
   v    &=& O_1 + r_1 e^{i (\pi/2+(\sgn\kappa_1)(\theta_1-\pi /2))} ,
\end{eqnarray*}
   one can calculate successively the parameters  
\begin{eqnarray*}
 \sgn\kappa_2 &=&  \sgn\sin((\sgn\kappa_1)\theta_1-(1/2+\alpha)\pi),\\ 
   O_2 &=& (\im v -\re v\tan
 ( (\sgn \kappa_1)\theta_1 + \frac{\sgn\kappa_2-1}{2}-\alpha\pi ) )\,i,\\
 \kappa_2 &=& \frac{\sgn\kappa_2}{|v-O_2|}, \\
   p_2 &=& O_2 + \frac{1}{\kappa_2}i, \\
 \theta_2 &=& (\sgn\kappa_2) \arg{\frac{p_2 - O_2}{v - O_2}} .
\end{eqnarray*}

If one knows only $O_2$ and $\theta_2$, there may be two
possible upper boundary edges for $D$, which are complementary arcs
of the circle centered at $O_2$ and passing through $v$. Thus
$\sgn\kappa_2$ is an essential datum. These formulas can also
be applied when $\kappa_1=0$, with the convention that
$\sgn0=0$.

We have shown that $D$ is determined by $p_1$, $\kappa_1$, $\theta_1$,
and $\alpha$.  Further, if $p_1$ and $\kappa_1$ are replaced by $1$
and $p_1\kappa_1$, the figure is rescaled.  Thus for purposes of
defining conformal equivalence classes we may normalize $p_1=1$ and
conclude that there is a three-dimensional space of conformally
distinct \scq s.

However, arbitrary values of the parameters given as above will not
necessarily define an \scq\ To begin with, $\kappa_1$ and $\theta_1$
must be such that the open right edge lies entirely in the right half
plane; for $\kappa_1>0$ this requirement is 
$p_1-r_1(1-\cos\theta_1)\ge0$, while for $\kappa_1<0$ any value of
$\theta_1$ in $(0,\pi]$ is possible.  Further, $\alpha$ must be chosen
so that $p_2$ does not lie in the lower half-plane. When 
$\kappa_1>0$, the upper edge leaving $v$ must not point simultaneously
into the right half-plane and the lower half-plane, so if
$\kappa_1>0$, the largest $\alpha\pi$ can be is $3\pi/2+\theta_1$, in
which case the upper edge of $D$ goes horizontally from $v$ to
$+\infty$ and continues from $-\infty$ to $-\cg{v}$.  When
$\kappa_1<0$, there are two cases. When $\alpha<1/2$, $\theta_1$ may
attain its maximum value of $\pi$, in which case $v$ coincides with
the vertex $\cg{v}$.  When $\alpha>1/2$ the largest possible
value of $\theta_1$ is $3\pi/2-\alpha\pi$, and again the upper edge of
$D$ goes horizontally through $\infty$.  We will loosely refer to
these \scq s as ``of extremal type'', when either
$v=\cg{v}\in\R$ or $\im v=\im p_2$.

\subsection{Schwarzian derivative and accessory parameters}
Let $\D=\{|z<1\}$ denote the unit disk.  Let $D$ be an \scq\ and let
$f\colon\D\to D$ be a conformal mapping such that $f(0)=0$, $f'(0)>0$.
The preimages of the vertices of $D$ are represented as 
$\pm e^{\pm it}$ where $0<t<\pi/2$.  We will write $\D_t$ when we wish
to emphasize the topological quadrilateral whose boundary is
$\partial\D$ with these four prevertices prescribed, in the context of
conformal mappings from the disk to another topological quadrilateral.

The Schwarzian derivative \cite{Neh} $\S_f$ of $f$ can be expressed as
$\S_f(z)=R_{\alpha,t,\lambda}(z)$ where
\begin{eqnarray}
  R_{\alpha,t,\lambda}(z) &=&   \label{eq:defR}
     4(1-\alpha^2) \, \psi_{0,t}(z)-2\lambda \, \psi_{1,t}(z), \\
  \psi_{0,t}(z) &=&   \nonumber
    \frac{  (\cos2 t)z^4  -2 z^2 + \cos2 t}
         {(z^4-2 (\cos 2 t) z^2+1)^2} , \\   
  \psi_{1,t}(z) &=& \frac{2\sin2 t}{z^4-2(\cos2 t)z^2+1}.\nonumber
\end{eqnarray}
(in the Appendix of \cite{Phoro} the justification for general
circle-arc polygons is worked out in detail; however
due to a typographic error $(1-\alpha)^2$ appears in place of
$(1-\alpha^2)$.) Observe that the parameters $\alpha$ and $\lambda$ in
$R_{\alpha,t,\lambda}(z)$ separate nicely in the formula
(\ref{eq:defR}); this fact will be important in the sequel.

A basic question is to investigate the set of combinations of
accessory parameters $(\alpha,t,\lambda)$ for which $f$ is a univalent
mapping.  The first parameter is the easiest to relate to the geometry
of the \scq\ $D$, being simply $1/\pi$ times the internal angle at the
vertices.  We will see that relationship between $t,\lambda$ and the
other geometric parameters (for instance $\kappa_1,\theta_1$) also
involves $\alpha$. It is convenient to do as much as possible with $t$
alone, without considering $\alpha$ and $\lambda$; this will be
accomplished in the next section.

\section{Mapping from a rectangle \label{sec:rectanglemap}}

In this section we carry out a change of variable, writing $f(z)=g(w)$
where $g\colon R\to D$ is a conformal mapping from a Euclidean
rectangle $R$ to an \scq\ $D$. Taking the domain to be $R$ instead of
$\D$ naturally changes the behavior of the mapping at the vertices.
In \cite{Pthesis} such mappings from a rectangle were studied for
$\alpha=0$ but without consideration of their relationship to the
corresponding maps from a disk.

We will fix the prevertices of $f$ choosing $0<t<\pi/2$. For $z\in\D$
define the complex elliptic integral
\begin{equation} \label{eq:defE}
 E(z) = \int_0^z \frac{dz}{ \sqrt{
         (z-e^{it})(z+e^{-it})(z+e^{it})(z-e^{-it})} } .
\end{equation}
Note that the values $\omega_1=E(1)$ and $\omega_2=E(i)$ are
respectively real and imaginary; it is well known that the image     
\begin{equation} \label{eq:rectangle}
 R = E(\D_t)=\{w\colon\ -\omega_1<\re w< \omega_1,\ -|\omega_2|<\im w<|\omega_2| \}
\end{equation}
is a rectangle with vertices at $\pm \omega_3$, $\pm\cg{\omega}_3$,  
where $\omega_3=\omega_1+\omega_2$.  We define 
\begin{equation} \label{eq:defg}
   g = f\circ E^{-1} ,
\end{equation}
the unique conformal mapping from $R$ to $D$ such that
$g(0)=0$, $g'(0)>0$.  By construction $g$ takes vertices of $R$ to 
vertices of $D$.

\subsection{The Schwarzian derivative and the $\wp$-function}
Our main task is to calculate the Schwarzian derivative of $f$, which
by the Chain Rule \cite{Neh} is equal to  
\begin{equation} \label{eq:Sgh}
  \S_f = (\S_g\circ E)E'^2 + \S_E.
\end{equation}
 When $f$ is continued analytically by reflection of $z$
across edges of $\D_t$ (i.e., maximal arcs of
$\partial\D-\{\mbox{prevertices}\}$), the image value $w=E(z)$ is
reflected along edges of $R$.  Even numbers of reflections return $z$
to its original value and effect translations of $w$. Since the
Schwarzian derivative annihilates M\"obius transformations, $\S_E$ is
a single valued function on the Riemann sphere and must therefore be a
rational function.

To describe $\S_g$ and $\S_f$ explicitly, let us introduce the auxiliary
function  
\begin{equation} \label{eq:defphi}
  \varphi(z) = \wp(E(z)+\omega_3)
\end{equation}
where $\wp=\wp_{\omega_1,\omega_2}$ is the Weierstrass
$\wp$-function \cite{WW} with primitive periods $2\omega_1$,
$2\omega_2$.  

As customary, write $e_i=\wp(\omega_i)$, $i=1,2,3$. These values are
real and satisfy
\begin{equation}  \label{eq:e2e1e3}
   e_2 <e_3<e_1.
\end{equation}
Let $Q$ denote the quarter-disk $Q=\{z\in\D\colon\ \pi<\arg
z<3\pi/2\}$, considered as a topological quadrilateral with vertices
$(-e^{it},-i,0,-1)$.  Since $\wp$ sends the subrectangle of $R$ with
vertices at $(0,\omega_1,\omega_3,\omega_2)$ to the lower half-plane,
it follows from (\ref{eq:defphi}) that $\varphi(Q)$ is also the lower
half-plane.  By applying the Reflection Principle one sees that
$\varphi$ must be a rational function of order $4$ and so we can write
it as
\[ \varphi(z) = a_0 + \frac{a_1}{z-e^{it}} + 
 \frac{a_2}{z+e^{-it}} + \frac{a_3}{z+e^{it}} + \frac{a_4}{z-e^{-it}} .
\] 
Clearly   $\varphi(-z)=\wp(E(-z)+\omega_3)=
\wp(-E(z)-\omega_3)=\wp(E(z)+\omega_3)=\varphi(z)$ and similarly
$\varphi(\cg{z})=\cg{\varphi(z)}$. From these symmetries
we find that $a_3=-a_1$, $a_4=-a_2$, $a_2=-\cg{a_1}$. Further,
\[ a_0 = \lim_{z\to\infty} \varphi(z)
       = \wp(E(\infty)+\omega_3) = \wp(\omega_3)=e_3
\]
since $E(\infty)$ is a period. We have so far
\[   \varphi(z) = e_3 + 4\frac{-z^2\re(a_1e^{it})+ \re(a_1e^{-it})}
                 {z^4-2(\cos2t)z^2+1} .
\]
Since $\varphi(0)=\wp(\omega_3)=e_3$ we have  
$\re(a_1e^{-it})=0$. From this it follows that
$\varphi''(0)=8\,\re(a_1e^{it})$. On the other hand, differentiating
the definition (\ref{eq:defphi}) and using the classical facts
$\wp'(\omega_3)=0$, $\wp''(\omega_3)=2(e_3-e_1)(e_3-e_2)$ together
with $E'(0)=1$, $E''(0)=0$ we arrive at 
\begin{equation} \label{eq:phi_1}
  \varphi(z) = e_3 - \frac{(e_1-e_3)(e_3-e_2)z^2}{z^4-2(\cos2t)z^2+1}.
\end{equation}
 
Since $\varphi$ is completely determined once $t$ is specified, we are
led to seek a formula for $\varphi(z)$ which does not involve  
$e_1,e_2,e_3$.  By the symmetry of $E$ with respect to reflections in
the coordinate axes, we know that $\varphi(-i)=e_1$,
$\varphi(-1)=e_2$, which reduces to
\begin{eqnarray*}
   e_1 &=& e_3 + 4 \sin^2 t ,\\
   e_2 &=& e_3 - 4 \cos^2 t .
\end{eqnarray*}
The resulting relation $e_1-e_2=4$, while perhaps surprising, is a
consequence of the normalization we have implicitly applied to the
periods of $\wp$.  From the well-known fact $e_1+e_2+e_3=0$ we
deduce that  
\begin{eqnarray} \label{eq:e1e2e3}
  e_1 &=&   2 - \frac{2}{3}\cos2t ,\\
  e_2 &=&  -2 - \frac{2}{3}\cos2t , \nonumber\\
  e_3 &=&  \frac{4}{3}\cos2t ,\nonumber 
\end{eqnarray}
and (\ref{eq:phi_1}) takes the definitive form
\begin{equation} \label{eq:phiz}
 \varphi(z)=  \frac{4}{3}\cos2t 
 -\frac{4  (\sin ^2 2t) z^2}{z^4-2 (\cos2t) z^2+1} .
\end{equation}

We observe that the half-period $\omega_1=K$ can be identified with
the complete elliptic integral in the terminology of Jacobian elliptic
functions \cite{WW}, and that $(e_3-e_2)/(e_1-e_2)=\cos^2t=k^2$ is the
modulus of the function $\mbox{sn}(z|k^2)$ which is related to $\wp$ by
\begin{equation} \label{eq:psn}
\wp(z)=e_2+(e_1-e_2)/\mbox{sn}^2(\sqrt{e_1-e_2}\,z)
\end{equation}
which in turn by (\ref{eq:e1e2e3}) is equal to
$(-2-(2/3)\cos2t)+4/\mbox{sn}^22z$.

Next we observe that $g$ extends by reflection along the edges of $R$
and is also symmetric under reflection in the coordinate axes. From
this it follows that when $\omega$ is a period of $\wp$, the value 
$g(\zeta+\omega)$ is related to $g(\zeta)$ by a M\"obius transformation.
Therefore $\S_g$ is doubly-periodic with periods $2\omega_1$ and
$2\omega_2$.  A simple local calculation verifies that $\S_g$ has
double poles at the vertices of $R$, and since $g$ converts the right
angles there into angles of $\alpha\pi$, we conclude that \ 
\begin{equation} \label{eq:Sg}
 \S_g(\zeta) = \frac{1-4\alpha^2}{2}\wp(\zeta+\omega_3) - 2\mu  
\end{equation}
for some constant $\mu\in\R$.  It is essential for us to see how $\mu$
depends on $t$ and on the geometry of the \scq\ $D$.  From  
(\ref{eq:defE})  we find that
\[ \S_E(z) = \frac{ 2(\cos2t)z^4+(\cos4t-5)z^2+2\cos2t}
       {(z^4-2(\cos2t)z^2+1)^2}.
\]
From this and (\ref{eq:Sgh}),(\ref{eq:defphi}), we obtain  
\begin{eqnarray}  \label{Sftemp}  
 \S_f &=&  (\frac{1-4\alpha^2}{2}\varphi-2\mu)E'^2  +\S_E   \\ 
      &=& \frac{(2/3)(1 - 4\alpha^2)\cos2t - 2\mu}{z^4-2(\cos2t)z^2+1}
     \nonumber\\ 
      &&  +\  \frac{ 2(\cos2t)z^4  + (\cos4t-2(1-4\alpha^2)\sin^22t-5)z^2
                    + 2\cos2t}
            { (z^4 - 2(\cos2t)z^2 + 1)^2 }\nonumber 
\end{eqnarray}
after some manipulation.
By comparison with (\ref{eq:defR}) we discover that
\begin{equation} \label{eq:mu}
 \mu = \frac{2(\alpha ^2-1)}{3} \cos2t + 2\lambda\sin2t .
\end{equation}
 Thus for fixed $\alpha$ and $t$, the correspondence
$\lambda\leftrightarrow\mu$ is linear, and $\lambda,\mu$ increase
simultaneously inasmuch as $0<t<\pi/2$. 

When one has generated a mapping $f$ via parameters $\alpha$, $t$,
$\lambda$, the corresponding parameters for $g$ in (\ref{eq:Sg}) are
obtained by keeping the same value of $\alpha$, using $\mu$ given by
(\ref{eq:mu}), and evaluating the integrals $E(1)$, $E(i)$ to
determine the defining periods for $\wp$. On the other hand, 
when working with (\ref{eq:Sg}) the data may be $\alpha$,
$\omega_2/\omega_1$ and $\mu$. Then one may proceed as follows.
 
We will say that a pair of half-periods $\omega_1\in\R$, $\omega_2\in i\R$
is \textit{normalized} when the corresponding values of the
$\wp$-function satisfy $e_1-e_2=4$. A given ratio
$\tau=\omega_2/\omega_1\in i\R$ is realized by a unique normalized
pair $(\omega_1,\omega_2)$.  To calculate it, one considers the
$\wp$-function with (non-normalized) primitive periods $1$, $\tau$
evaluated at the half-periods and calculates
$r=\sqrt{\wp(1/2)-\wp(\tau/2)}/4$, and then takes $\omega_1=r$,
$\omega_2=r\tau$.  This may be done easily via Jacobian theta
functions, using the formulas
\begin{eqnarray*}
\wp(1/2) &=& \phantom{-} \frac{\pi^2}{3}( \vartheta_3(0)^4 + \vartheta_4(0)^4),\\
 \wp(\tau/2) &=& -\frac{\pi^2}{3}( \vartheta_3(0)^4 + \vartheta_2(0)^4),
\end{eqnarray*}
in the notation of \cite{WW}, with the elliptic modulus 
$q=e^{i\pi\tau}$.

We observe that with normalized periods, as $\tau$ runs through the
imaginary axis, the triple $(e_1,e_2,e_3)=(2-e_3/2,\ -2-e_3/2,\ e_3)$
traces the straight segment in $\R^3$ from $(4/3, -8/3, 4/3)$ to
$(8/3,-4/3,-4/3)$.

Our point of view for the conformal mapping problem is that although
the Schwarzian derivative (\ref{eq:Sg}) is a transcendental function,
in some respects it is simpler to study than the rational function
(\ref{eq:defR}).

\subsection{The zeros of $\wp$ expressed as elliptic integrals
\label{subsec:zeros}}

We digress to observe that this setting sheds some light on the
location of the zeros of $\wp$.  In \cite{DI} a formula for the zeros
of $\wp$ is given in terms of generalized hypergeometric
functions. The starting point for the derivation of this formula is a
formula in \cite{EZ} which is an improper integral involving the well
known modular forms $E_6(\tau)$ and $\Delta(\tau)$. The latter formula
may in turn be derived from the elliptic integral
\begin{equation} \label{eq:z0a}
 z_0(\tau)=\pm \int_0^\infty {dx \over \sqrt{4x^3-g_2x-g_3}}\nonumber
\end{equation}
for the zeros $z_0(\tau)$ of $\wp_\tau$, where $4x^3-g_2x-g_3=
4(x-e_1)(x-e_2)(x-e_3)$. 

We make the observation (not mentioned in \cite{DI} and \cite{EZ})
that another formula for $z_0$ in terms of an elliptic integral is
\begin{equation} \label{eq:z0b}
  z_0(\tau) = \pm {1\over \sqrt{e_1-e_2}}  \int_0^{2\sqrt{-1/e_2}}
   \frac{dv}
     {\sqrt{\vphantom{({a\over b})^2} 1-v^2}
      \sqrt{1-( \frac{e_3-e_2}{e_1-e_2})v^2}}  .
\end{equation}
This formula is an immediate consequence of the formula
(\ref{eq:psn}).  We derive yet another formula for $z_0$ as an
elliptic integral: it follows from (\ref{eq:defphi}) that if
$\pi/4<t<\pi/2$ then the assertion $\wp(z_0)=0$ can be expressed as
\mbox{$\varphi(-iy_0)=0$}, where  
\[ z_0=E(-iy_0)+\omega_1+\omega_2 ,
\] 
for some real number $y_0$  in the interval $0<y_0<1$. Thus
by (\ref{eq:phiz}),
\[
 {4\over 3}\cos2t+{4(\sin^22t)y_0^2\over y_0^4+2(\cos2t)y_0^2+1}=0 ,
\]
which implies
\begin{equation}  \label{eq:y02}
 y_0^2 = \frac{ -3+\cos^22t+(\sin2t)\sqrt{9-\cos^22t} }{2\cos2t} .
\end{equation}
We therefore have
\begin{eqnarray}
 z_0 &=& \omega_1+\omega_2 -i\int_0^{y_0} 
   \frac{dx}{\sqrt{x^4+2(\cos2t)x^2+1}}  \nonumber \\
    &=&  \omega_1+\omega_2-{i\over 2} \int_0^{{2y_0 /(1+y_0^2)}}
     \frac{dv}{ \sqrt{\vphantom{({a\over b})^2} 1-v^2} 
     \sqrt{1-(\sin^2t)v^2} }  . \label{eq:z0c}
\end{eqnarray}
By (\ref{eq:y02}) we find that
\[ {4y_0^2 \over (1+y_0^2)^2 } = 
     \frac{4\cos2t}{\cos^22t +2\cos2t-3} =
     \frac{16e_3 }{ (e_3+4)(3e_3-4) }  .
\]
Consequently, (\ref{eq:z0c}) becomes
\begin{equation}\label{eq:z0d}
  z_0 = \omega_1+\omega_2 - \frac{i}{2} 
   \int_0^{4\sqrt{ e_3 / ( (e_3+4)(3e_3-4) )}}
   \frac{dv}{\sqrt{\vphantom{({a\over b})^2} 1-v^2} 
             \sqrt{1-({e_1-e_3 \over 4})v^2 }}\, . 
\end{equation}
By analytic continuation, the formulas (\ref{eq:z0b}) and
(\ref{eq:z0d}) are valid in general; i.e.,
$\wp(z_0(\tau)\,|\,\tau)=0$ is valid for $z_0(\tau)$ computed using
(\ref{eq:z0b}) or (\ref{eq:z0d}) with $\im\tau>0$, since
$e_1(\tau),e_2(\tau),e_3(\tau)$ and $z_0(\tau)$ are analytic functions
of $\tau$.

\section{Calculation of the image quadrilateral as a function
of the accessory parameter \label{sec:calc}}

In this section we will carry further the method initiated in \cite{KP2}
by which geometric characteristics of the image \scq\ are represented
as power series in the accessory parameter $\lambda$ of (\ref{eq:defR})  
or $\mu$ of (\ref{eq:Sg}). This is a straightforward application of the
recent Spectral Parameter Power Series (SPPS) method for solution of 
Sturm-Liouville differential equations \cite{KP1}, which we restate here
briefly to make the presentation self-contained.

\subsection{SPPS method}
We will work on the real interval $0\le z\le1$.  One defines the
sequence $I_n$ of \textit{iterated integrals} generated by an
arbitrary pair of functions $(q_0,q_1)$ recursively by setting $I_0=1$
identically and for $n\ge1$,  
\begin{equation}
  I_n(z) = \int_0^z I_{n-1}(\zeta)\,q_{n-1}(\zeta)\,d\zeta
\end{equation}
where $q_{n+2j}=q_n$ for $j=1,2,\dots$

\begin{proposition}{\rm \cite{KP1}} \label{prop:powersX}
 Let\/ $\psi_0$ and $\psi_1$ be given, and suppose that\/ $y_\infty$
is a nonvanishing solution of
\[ y_\infty''+\psi_0\,y_\infty=\lambda_\infty\psi_1\,y_\infty
\] on the interval $[0,1]$.
Choose\/ $q_0=1/y_\infty^2$, $q_1=\psi_1\,y_\infty^2$ and define\/
$\X{n}$, $\Xt{n}$ to be the sequences of iterated integrals generated
by\/ $(q_0,q_1)$ and by\/ $(q_1,q_0)$ respectively.  Then for each\/
$\lambda\in\C$ the functions
\begin{eqnarray*}
  y_1 = y_\infty \sum_{k=0}^\infty(\lambda-\lambda_\infty)^k
            \Xt{2k},\\ 
  y_2 = y_\infty \sum_{k=0}^\infty(\lambda-\lambda_\infty)^k
              \X{2k+1}\\
\end{eqnarray*}
are linearly independent solutions of the equation
\begin{equation}  \label{eq:y''}
 y''+\psi_0y=\lambda\psi_1y 
\end{equation}
on  $[0,1]$.  Further, the series for $y_1$ and $y_2$ 
converge uniformly on $[0,1]$ for every $\lambda$.
\end{proposition}

Thus two linearly independent solutions of the second order ordinary
differential equation (\ref{eq:y''}), which depend on the parameter
$\lambda$, are expressed as explicitly calculable analytic
functions of $\lambda$.

\subsection{$p_2$ and $\kappa_1$  \label{subsec:p2kappa1}}

We now concentrate on the geometric parameters $p_2$ and $\kappa_1$
described in section \ref{sec:quads}, where the \scq\ $D$ is
normalized by $p_1=1$.  Suppose that $\alpha,t,\lambda$ are given.
The formulation we give here for general $\alpha$ is analogous to the
case $\alpha=1/2$ as presented in \cite{KP2}.  Let $f$ be the solution
of $\S_f=R_{\alpha,t,\lambda}$ in $\D$ (recall (\ref{eq:defR}))
normalized by
\begin{equation} \label{eq:f0}
 f(0)=0,\quad f'(0)=1,\quad f''(0)=0,
\end{equation}
and suppose that $f(\D)$ is an \scq\ similar to $D$, specifically
$f(\D)= w_1D$ where $w_1>0$. Further, as is common, take $y_1$ to be a
solution of (\ref{eq:y''}) in $\D$ with the normalization $y_1(0)=1$,
$y_1'(0)=0$, where we use $\psi_{0,t}$, $\psi_{1,t}$ in place of
$\psi_0$, $\psi_1$ in Proposition \ref{prop:powersX}. Then one may
calculate the right edge midpoint $w_1=f(1)$ via either of the
expressions in  
\begin{equation} \label{eq:fintegral}
 f(z) = \int_0^z \frac{1}{y_1(u)^2} \,du 
      = \frac{y_2(z)}{y_1(z)},
\end{equation}
where $y_2$ is the solution of (\ref{eq:y''}) normalized by
$y_2(0)=0$, $y_2'(0)=1$.  Similarly, the upper midpoint $w_2=f(i)$ can
be calculated by integrating along the segment from $0$ to $i$, or
more easily, by using $R_{\alpha,\pi/2-t,-\lambda}$ instead of
$R_{\alpha,t,\lambda}$ and integrating from $0$ to $1$.  We will refer
to the corresponding solutions of this latter integration as $y_1^*$,
$y_2^*$. Now midpoint of the upper edge is seen to be
\begin{equation}\label{eq:p2}
  p_2 = \frac{w_2}{w_1} = \frac{y_1(1)y_2^*(1)}{y_1^*(1)y_2(1)} .
\end{equation}

It was observed for example in \cite{Br2,KP2} that the curvature of
the right edge of $f(\D)=w_1D$ is equal to $y_1(1)^2 -
2y_1(1)y_1'(1)$.  This can be verified by writing out the relationship
$|f(z)-O_1||f(z^*)-O_1|=1/(\mbox{curvature})^2$ which relates points
$z,z^*$ symmetric with respect to the unit circle, and examining the
limit as $z\to 1$, $z^*\to 1$. The result holds for zero curvature as
well although the proof must be modified for that case.  The curvature
of the right edge of $D$ is $w_1$ times the value for $w_1D$,
\begin{equation} \label{eq:kappa}
  \kappa_1 = y_1(1)y_2(1)  - 2y_2(1)y_1'(1) .
\end{equation}

Let $\alpha,t$ be fixed.  Then how do $p_2$, $\kappa_1$ depend upon
$\lambda$? Suppose we know a parameter value $\lambda_\infty$ such
that the conformal mapping $f_{\alpha,t,\lambda_\infty}$ whose
Schwarzian derivative is $R_{\alpha,t,\lambda_\infty}$ is
univalent. (The question of finding such $\lambda_\infty$ will be
discussed in section \ref{sec:domuniv}). Then it follows that the
solution $y_\infty$ of (\ref{eq:y''}) with $\lambda=\lambda_\infty$
and with the normalization $y_\infty(0)=1$, $y_\infty'(0)=0$ is
nonvanishing in $\D$, and hence in $[0,1]$.  This $y_\infty$ may be
used to define the power series $y_1[\lambda]$, $y_2[\lambda]$ of
Proposition \ref{prop:powersX}.  The notation $[\lambda]$ refers here
to evaluation at $z=1$ of the function whose parameter is $\lambda$.
Then (\ref{eq:p2}) and (\ref{eq:kappa}) give us
\begin{eqnarray}
   p_2(\lambda) &=&  \label{eq:p2lambda}
      \frac{y_1[\lambda]y_2^*[\lambda]}{y_1^*[\lambda]y_2[\lambda]} \\
   \kappa_1(\lambda) &=& 
     y_1[\lambda]y_2[\lambda]  - 2y_2[\lambda]y_1'[\lambda]
 \label{eq:kappa1lambda}
\end{eqnarray}
For evaluation of $p_2(\lambda)$ at a single $\lambda$
it one may prefer to evaluate the four power series on the
right side of (\ref{eq:p2lambda}) separately rather than to work
out a single series.  For solving equations of the form
$p_2(\lambda)=c$, one may rewrite this equation as
$y_1[\lambda]y_2^*[\lambda]-cy_1^*[\lambda]y_2[\lambda]=0$, 
thus avoiding division of power series.

Now we consider the map $g\colon R\to D$ given by (\ref{eq:defg}),  
defined in the rectangle $R$. Let $\eta_1$ and $\eta_2$ be the solutions of
\begin{equation}
   \eta''(\zeta) + \left( (\frac{1}{4}-\alpha^2)\wp(\zeta+e_3) - 
   \mu \right) \eta(\zeta) = 0 
\end{equation}
with normalizations $(\eta(0),\eta'(0))=(1,0)$ and $(0,1)$
respectively, so 
\[ g =  \frac{\eta_2}{\eta_1}
\]
because $g'(0)=f'(0)=1$ (recall that $E'(0)=1$ in section
\ref{sec:rectanglemap}). We have $g(\omega_1)=w_1$, $g(\omega_2)=w_2$
so it is easy to express $p_2$.  As to the curvature, 
we note that for $\zeta\in R$, the symmetry is
\[  |g(\zeta)-O_1| |g(2\omega_1-\zeta)-O_1| = 
     \frac{1}{\mbox{curvature}^2}
\]
from which we find the curvature of the
right edge of $w_1D$ to be $2\eta_1(\omega_1) \eta_1'(\omega_1)$,
which is a simpler expression than for the mapping $f$ defined in the
disk. Thus
\begin{equation} \label{eq:kappaeta}
 \kappa_1 =  \frac{\eta_2(\omega_1)}{\eta_1(\omega_1)}
         (2\eta_1(\omega_1) \eta_1'(\omega_1)) 
    = 2\eta_2(\omega_1) \eta_1'(\omega_1).
\end{equation}
Suppose that $\alpha$ and the normalized pair of half-periods  
$\omega_1,\omega_2$ (or their ratio $\tau$) are given.  Analogously to
the case of mapping from $\D$, if we have a parameter $\mu_\infty$ for
which $g$ is univalent, we can form power series in $\mu$  
\begin{eqnarray}
   p_2(\mu) &=&  \label{eq:p2mu}
      \frac{\eta_1[\mu]\eta_2^*[\mu]}{\eta_1^*[\mu]\eta_2[\mu]}, \\
   \kappa_1(\mu) &=&  \label{eq:kappamu}
   2\eta_2[\mu] \eta_1'[\mu].
\end{eqnarray}
centered at $\mu_\infty$, where $\eta_1^*[\mu],\eta_2^*[\mu]$ refer to
the values obtained by integrating the differential equation
vertically from $0$ to the endpoint $\omega_2$.

From the above we are able to calculate $\kappa_1$, $p_2$ numerically
given the accessory parameters of the Schwarzian derivative, both for
maps of $\D$ and of $R$.  Further, if we assume that $\alpha,t$ or
$\alpha,\tau$ are fixed and one of the values $\kappa_1$ or $p_2$ is
specified, then we can approximate the parameter $\lambda$ or $\mu$
arbitrarily closely by truncating the corresponding power series to a
polynomial of sufficiently high degree and calculating an appropriate
real zero.  Finally, when the full geometry of $D$ is specified via
$\alpha$, $\kappa_1$, and $p_2$, we can calculate $(t,\lambda)$ or
$(\tau,\mu)$ by an iterative procedure in which a guessed value for
the modulus $t$ or $\tau$ is corrected by solving for $\lambda$ or
$\mu$ to achieve either one of $\kappa_1$ or $p_2$, and then adjusting
to make the other value come out correct as well.  There is no
essential difference to the way this was done in \cite{KP2} for
$\alpha=1/2$, so we will not give further details here.

\section{Parameters for univalence\label{sec:domuniv}}

\begin{figure}[t!]
\pic{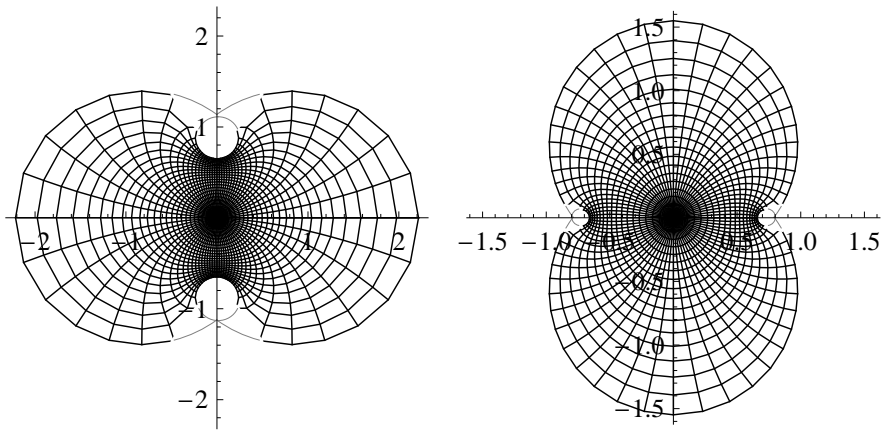}{6.5,0}{4.5}{scale=1.1}
\caption{\small Extremal domains for $\alpha=0.2$, $t=\pi/6$, with
$\lambda_{\rm min}=-0.479608$ (left) and   $\lambda_{\rm max}=1.30611$
(right). Values near the vertices were plotted by extrapolation.
\label{fig:lamminmax1}}
\end{figure}
 
\begin{figure}[b!]
\pic{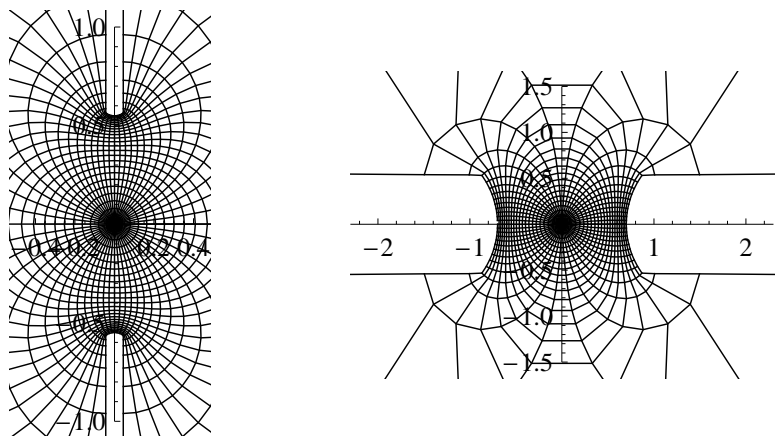}{6.5,0}{5.5}{scale=1.1}
\caption{\small Extremal domains for $\alpha=1.3$, $t=3\pi/8$, with
$\lambda_{\rm min}=-1.28089$  (left) and   $\lambda_{\rm max}= 1.84854$ 
(right).
\label{fig:lamminmax2}}
\end{figure}

\subsection{Interval of univalence for fixed $\alpha$,$t$ 
\label{subsec:intervuniv}}
 
Following the line of investigation of \cite{Br1,Br2,KP2}, let us
assume that $\alpha$,$t$ are given.  We investigate the interval
$\lambda_{\rm min}\le\lambda\le\lambda_{\rm max}$ of values of
$\lambda$ for which the map $f=f_{\alpha,t,\lambda}$ is univalent.  A
picture of the set of univalence pairs $(t,\lambda)$ corresponding to
$\alpha=1/2$ was given in \cite{KP2}.

By the symmetry of \scq s, a univalent $f$ cannot have a pole in
$\D$.  As an aid to the discussion we introduce the interval
$\lambda^*_{\rm min}<\lambda<\lambda^*_{\rm max}$ of values of $\lambda$
for which $f$ has no poles in $\D$. Thus
\begin{equation} \label{lambdas}
  \lambda^*_{\rm min} \le \lambda_{\rm min} <\lambda_{\rm max}  \le
   \lambda^*_{\rm max}.
\end{equation}
Recalling the discussion of domains of extremal type in section
\ref{subsec:descrip}, we first suppose that $\alpha\le1/2$.  Then
$\lambda=\lambda_{\rm max}$ produces an extremal domain when the
vertices $v=f(e^{it})$ and $f(e^{-it})$ coincide on the positive real
axis.  When this happens, $v=w_1+2r_1$ where as previously, $w_1=f(1)$
and $\kappa_1=-1/r_1<0$. From elementary geometry it is seen that
the center of the upper edge of the \scq\ $D$ is $O_2=w_2-ir_2$ 
and that $\arg(v-O_2)=-\alpha\pi$.  From this it follows that
\[   \sin\pi\alpha = \frac{|w_2|-r_2}{r_2}
\]
which says
\[   |w_2|\kappa_2- \sin\alpha\pi - 1=0 
\]
since $\kappa_2=1/r_2>0$.  By section \ref{sec:calc} we know how to
express the left hand side of this equation as a power series in

$\lambda$, which we may solve numerically for $\lambda_{\rm max}$.
(The art of locating the zeroes of a polynomial with real coefficients
is quite well developed and we will not discuss it here.)  Similarly,
the value $\lambda=\lambda_{\rm min}$ for which the upper and lower
edges of $D$ degenerate into full circles may be calculated by
applying the same procedure using $\pi/2-t$, $-\lambda$ in place of
$t$,$\lambda$. Once $\lambda_{\rm max}$ (or $\lambda_{\rm min}$) is
found, it is easy to calculate the conformal mapping. Examples are
given in Figure \ref{fig:lamminmax1} ($\alpha<1/2$) and Figure
\ref{fig:lamminmax2}  ($\alpha>1/2$).

The remaining case $\alpha>1/2$ is actually simpler.  The upper
and lower edges of $D$ degenerate into horizontal rays meeting at
$\infty$ precisely when $f(i)=\infty$, i.e., we have $\lambda^*_{\rm
  max}=\lambda_{\rm max}$.  The degeneration of the right and left
edges, $f(1)=\infty$, corresponds to $\lambda^*_{\rm min}=\lambda_{\rm
  min}$.

\subsection{Estimates for $\lambda^*_{\rm min}$, $\lambda^*_{\rm max}$.}

There are certain advantages to working with the rectangle map $g$
because the transcendental function $\S_g$ in many ways behaves more
nicely than the rational function $\S_f$. Thus we will investigate the
values $\mu^*_{\rm min},\mu^*_{\rm max}$ which correspond to
$\lambda^*_{\rm min},\lambda^*_{\rm max}$ by (\ref{eq:mu}).  Since
$\wp$ is real on the coordinate axes as well as on the boundary of the
rectangle  with vertices $\pm \omega_3$, $\pm\cg{\omega}_3$,
and since $\wp'$ only vanishes at half-periods which are not periods,
it follows from (\ref{eq:Sg}) that the real function $\S_g$ is
one-to-one on the horizontal segment from $0$ to $\omega_1$ and also
on the vertical segment from $0$ to $\omega_2$.  This can be made more
precise noting that $\wp$ decreases from $e_3$ to $e_2$ on
$[\omega_3,\omega_3+\omega_1]$ while it increases from $e_3$ to $e_1$
on $[\omega_3,\omega_3+\omega_2 ]$. In contrast, it appears to be much
more difficult to determine bounds for $S_f$ on the coordinate axes in
$\D$.

\begin{lemma} \label{lem:sturm}  
Let $\eta''(\zeta)+a(\zeta)\eta(\zeta)=0$ on a real interval  
$[\zeta_0,\zeta_0+L]$, and let $\eta(\zeta_0)=1$,
$\eta'(\zeta_0)=0$. (i) Suppose $a(\zeta)\le0$ for all 
$\zeta\in[\zeta_0,\zeta_0+L]$. Then $\eta$ never vanishes. (ii)
Suppose $a(\zeta)\ge(\pi/(2L))^2$ for all $\zeta$. Then $\eta$
vanishes at least once in the interval $[\zeta_0,\zeta_0+L]$.
\end{lemma}

\proof
The proof is immediate from the Sturm comparison theorem \cite{Har}.
For (i), note that the similarly normalized solution of
$\eta''(\zeta)+0\cdot\eta(\zeta)=0$ is $\eta\equiv1$, which never vanishes.
For (ii), use comparison with 
$\eta''(\zeta)+(\pi/(2L))^2\cdot\eta(\zeta)=0$, whose solution
$\cos((\pi/(2L))(\zeta-\zeta_0))$ vanishes at $\zeta_0+L$. \qed

\begin{proposition} \label{prop:mubounds} Let the accessory parameters 
$\alpha$ and $t$ (or equivalently $\tau$) be fixed.  (i) If
  $\alpha\le1/2$, then  
\[   (\frac{1}{4}-\alpha^2)e_2-\left(\frac{\pi}{2\omega_1}\right)^2
      \le \mu^*_{\rm min}, 
\] \[   \mu^*_{\rm max}
   \ \le \ (\frac{1}{4}-\alpha^2)e_1 + 
     \left(\frac{\pi}{2|\omega_2|}\right)^2  .
\]
(ii) If $\alpha\ge1/2$, then  
\[   (\frac{1}{4}-\alpha^2)e_3-\left(\frac{\pi}{2\omega_1}\right)^2
      \le \mu^*_{\rm min},   
\] \[    \mu^*_{\rm max}
   \ \le \ (\frac{1}{4}-\alpha^2)e_3 + 
     \left(\frac{\pi}{2|\omega_2|}\right)^2  .
\]
\end{proposition}

\proof From (\ref{eq:Sg}) we are led to look at $g=\eta_2/\eta_1$ with
$\eta_1$ the solution of   
\begin{equation} \label{eq:eta''}
 \eta''(\zeta) + ((\frac{1}{4}-\alpha^2)\wp(\zeta+\omega_3)-\mu)\eta(\zeta)=0
\end{equation}
on the interval $0\le\zeta\le\omega_1$, normalized as in Lemma
\ref{lem:sturm}. First suppose that
$\alpha\le1/2$. Then by (\ref{eq:e2e1e3}) the coefficient of $\eta$ in
(\ref{eq:eta''}) satisfies the bounds
\[   (\frac{1}{4}-\alpha^2)e_2-\mu \le 
     (\frac{1}{4}-\alpha^2)\wp(\zeta+\omega_3)-\mu   \le 
     (\frac{1}{4}-\alpha^2)e_3-\mu,
\]
so $\eta_1$ must vanish at some point if
$(\frac{1}{4}-\alpha^2)e_2-\mu> (\pi/(2\omega_1))^2 $, which tells us
that $\mu\le\mu^*_{\rm min}$ for such values of $\mu$.  Thus the first
inequality of (i) is established.  Further, we note that according to
the Lemma, $\eta_1$ will never vanish if
$\mu\ge(\frac{1}{4}-\alpha^2)e_3$.

To examine the behavior in the imaginary direction it is helpful to
introduce $\hat\eta(u)=\eta(\zeta)$ where $\zeta=iu$ for $0\le
u\le|\omega_2|$. The differential equation is
\begin{equation} \label{eq:etahat''}
 \hat\eta''(u) + 
      (\mu-(\frac{1}{4}-\alpha^2)\wp(iu+\omega_3))\,\hat\eta(u)=0,
\end{equation} 
so the normalized solution $\hat\eta_1$ must vanish   if
$\mu-(\frac{1}{4}-\alpha^2)e_1> (\pi/(2|\omega_2|))^2$, which tells us
that $\mu\ge\mu^*_{\rm max}$ for such values of $\mu$.  Thus the final
inequality of (i) is  established. Further, we note that  
$\hat\eta_1$ will never  vanish if
$\mu\le(\frac{1}{4}-\alpha^2)e_3$.

Part (i) has been proved.  The proof of (ii) is similar; the
main difference is that  
the inequalities derived from (\ref{eq:e2e1e3}) are reversed when
$1/4-\alpha^2$ is negative.  \qed

 From two facts which we noted in the course of the proof, $\eta$ and
 $\hat\eta$ do not both vanish for a common value of $\mu$; i.e., $g$
 does not simultaneously have poles on both coordinate axes. The
 relationship of the value $(\frac{1}{4}-\alpha^2)e_3$ to $\mu^*_{\rm
   min}$, $\mu^*_{\rm max}$, and univalence appears to merit further
 investigation.
 
From (\ref{lambdas}) we have immediately the following bounds on
$\mu_{\rm min}$, $\mu_{\rm max}$.
\begin{coro}
 For $\alpha\le1/2$, 
\[  (\frac{1}{4}-\alpha^2)e_2-\left(\frac{\pi}{2\omega_1}\right)^2
      \le \mu_{\rm min}, \quad
  \DS\mu_{\rm max}  \ \le \ (\frac{1}{4}-\alpha^2)e_1 + 
     \left(\frac{\pi}{2|\omega_2|}\right)^2;
\] 
whereas for $\alpha\ge1/2$,
\[ (\frac{1}{4}-\alpha^2)e_3-\left(\frac{\pi}{2\omega_1}\right)^2
      \le \mu_{\rm min}, \quad
  \DS\mu_{\rm max}  \ \le \ (\frac{1}{4}-\alpha^2)e_3 + 
     \left(\frac{\pi}{2|\omega_2|}\right)^2.
\]
\end{coro}

It follows from (\ref{eq:e2e1e3}) that  
$(1/4-\alpha^2)e_2-(\pi/(2\omega_1))^2 \le \mu_{\rm min}$ and
$\mu_{\rm max}\le(1/4-\alpha^2)e_1 + (\pi/(2|\omega_2|))^2$ for all
$\alpha$. The bounds can easily be converted to bounds for
$\lambda_{\rm min}$, $\lambda_{\rm max}$ by means of (\ref{eq:mu}).

\begin{figure}[!t]
 \pic{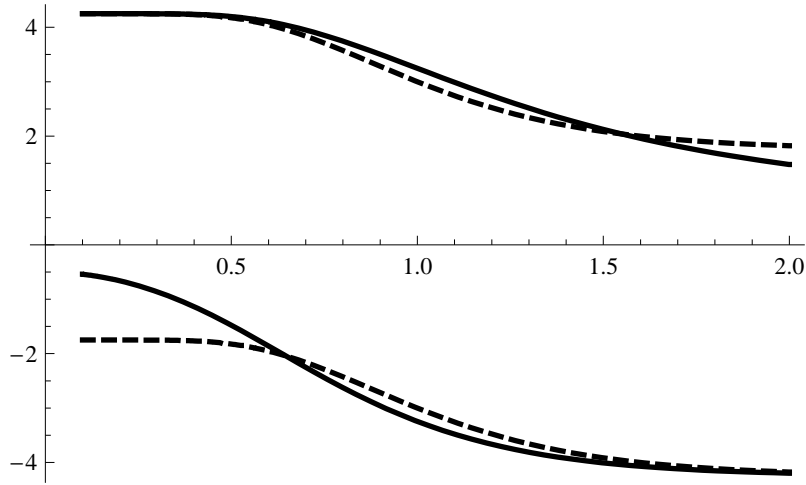}{7,.2}{6.1}{}
\caption{\small Comparison of upper and lower bounds obtained for
  $\mu_{\rm min}$, $\mu_{\rm max}$ for $\alpha=0.25$ and
  $0.1\le\tau/i\le2.0$, given by Proposition \ref{prop:mubounds}
  (solid lines) and Proposition \ref{prop:nehari} (dashed
  lines). \label{fig:mustarnehari}}
\end{figure}

 Bounds for
$\lambda_{\rm min}$, $\lambda_{\rm max}$ may also be derived in
another way, by using Nehari's theorem \cite{Neh1} which affirms that
$f$ holomorphic in $\D$ is not univalent when
\begin{equation}
   \sup_{z\in\D}(1-|z|^2)|\S_f(z)| > 6.
\end{equation}
Thus a necessary condition for univalence is $|\S_f(0)|\le6$, but
by (\ref{Sftemp}) we have  $\S_f(0)=2(1-\alpha^2)e_3-2\mu$, so we arrive
at the following result. 
\begin{proposition} \label{prop:nehari}
$(1-\alpha^2)e_3-3 \le\mu_{\rm min}$ and
$\mu_{\rm max} \le (1-\alpha^2)e_3+3$.
\end{proposition}
Numerically we can find values of $\alpha,t$ to make the bounds of
either one of Propositions \ref{prop:mubounds} and \ref{prop:nehari}
better than the other (an example is shown in Figure
\ref{fig:mustarnehari}).  None of the above arguments can be expected
to provide sufficient conditions for univalence, since our reasoning
only samples a small subset of the domain of the mapping.
 
 A natural candidate for a value of $\lambda$ within the range of
univalence is the average of the upper and lower bounds.
Numerical experiments were made with the average of the bounds given
by Propositions \ref{prop:mubounds} and \ref{prop:nehari}, as well as
averaging the combined bounds (the maximum of the lower bounds and the
minimum of the upper bounds).  All of these appear to provide
univalent mappings over a wide range of values of $\alpha$ and $t$,
the best results coming from the Nehari bounds of Proposition
\ref{prop:nehari}.

Once one has succeeded in locating a value $\lambda_\infty$ (or
$\mu_\infty$) within the range of univalence, i.e., for which $y_1$ is
a nonvanishing solution of the second order linear differential
equation, the SPPS method enables calculation of $\kappa_1$,$w_1$,
etc., to considerable accuracy.  (Actually, for the SPPS method to
work it is not necessary to be within the domain of univalence, but
only for the mapping to have no poles on the axis on which one is
integrating.)  Thus it is a simple matter to determine precise values
of $\lambda_{\rm min}$, $\lambda_{\rm max}$ numerically, as the
nearest zeros of an appropriate power series (approximated by a
polynomial) centered at $\lambda_\infty$.  This was done in \cite{KP2}
for $\alpha=1/2$ and there is no essential difference for general
$\alpha$.

\begin{figure}[!h]
\begin{center}
 \pic{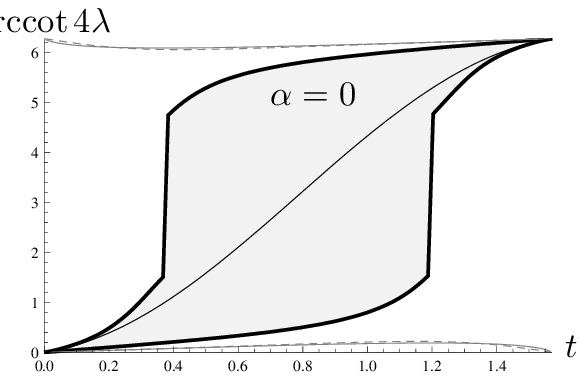}{-1,-2.2}{1.3}{} 
 \pic{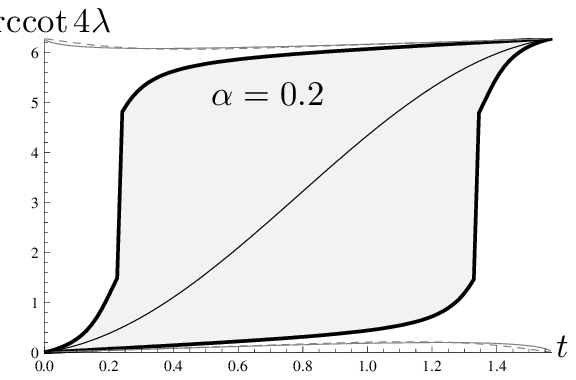}{5.2,-2.2}{0}{}\\
 \pic{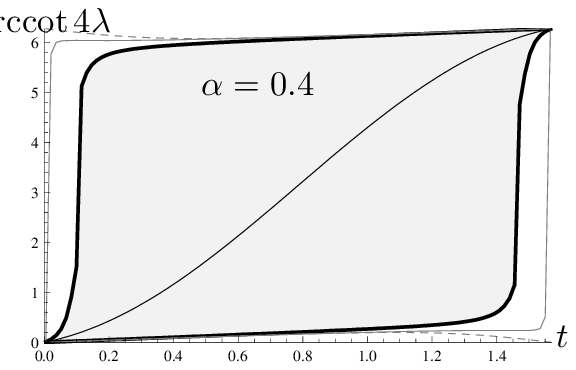}{-1,-1.5}{4.3}{} 
 \pic{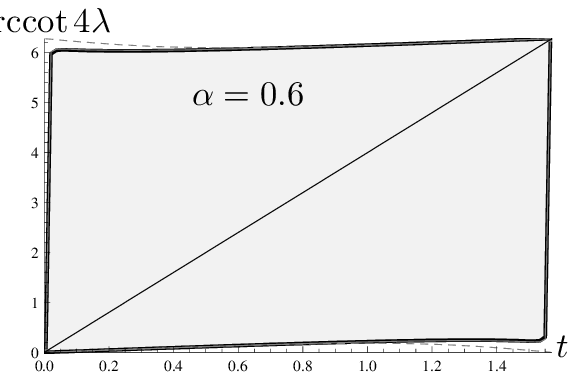}{5.2,-1.5}{0}{} \\
 \pic{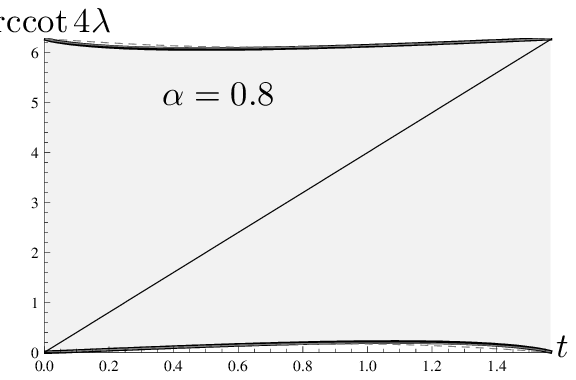}{-1,-.6}{4.3}{} 
 \pic{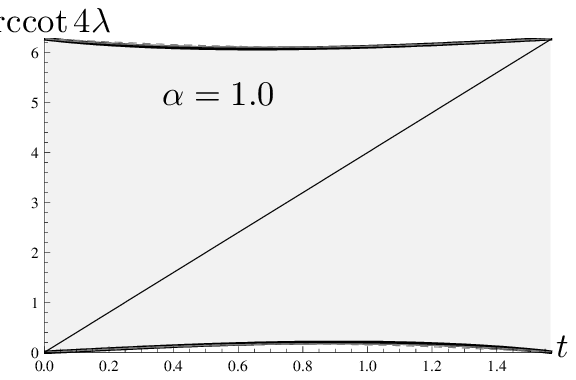}{5.2,-.6}{0}{}  
\end{center}
\caption{\small Values of $\mbox{arccot}\,4\lambda$ for which the mapping
  $f_{\alpha,t,\lambda}$ ($0<t<\pi/2$) is univalent in the unit
  disk. Bounds derived from $\mu^*_{\rm min}$,$\mu^*_{\rm max}$ are
  shown with a solid gray line; bounds derived from Nehari's theorem
  with a dashed line. \label{fig:univalence}}
\end{figure}

In Figure \ref{fig:univalence} we show the $\lambda$-domains of
univalence for $\alpha=0.0$, $0.2$, $0.4$, $0.6$, $0.8$, and $1.0$,
together with the bounds described above. (The pictures for  
$\alpha>1$
turn out similar to $\alpha=1$)  
Following the presentation in \cite{KP2}, we have applied the
transformation $\mbox{arccot}\,4\lambda$ to make the domain finite.
For $\alpha>1/2$ the bounds become difficult to distinguish from the
contours of the domains.  Domains for $\alpha>1$ are similar and are
not shown. The median curves dividing the these domains of univalence
in symmetric halves were devised as follows: The slopes of the
boundary curves at the corners $(0,0)$ and $(\pi/2,2\pi)$ were
estimated numerically based on the calculated values of $\lambda_{\rm
  min}$ or $\lambda_{\rm max}$; write $\beta$ for the slope at
$(0,0)$. Then the median curve is given by the formula
\[ \pi + \beta (t - \frac{\pi}{4}) - \pi (1 - \frac{\beta}{4}) \cos2t .
\]
This appears nicer than a simple average of the upper and lower
contours, and generalizes the dividing curve for $\alpha=1/2$   
which is a straight segment as shown in \cite{KP2}.  One may use
this median curve as starting values  $\lambda_\infty$ for the
SPPS method.
 
\section{Universal cover of a doubly-punctured disk\label{sec:dpd}} 

As an application of the calculation of the accessory parameters of
mappings of \scq s, we consider the triply-connected plane domain
$G_a=\D-\{-a,a\}$ where the value $0<a<1$ is fixed. The universal
covering map $H\colon\D\to G_a$ has been studied in \cite{Hem1,Hem2,SH}.
We can assume $H(0)=0$ and $H'(0)>0$. The first question is to
calculate or estimate $H'(0)$. Another is to calculate the
generators of the covering group of $G_a$ corresponding to the
covering $H$. We will use many basic facts about Fuchsian groups and  
covering spaces of Riemann surfaces, for which background information 
may be found in \cite{Mas}.

\begin{figure}[!b]
 \pic{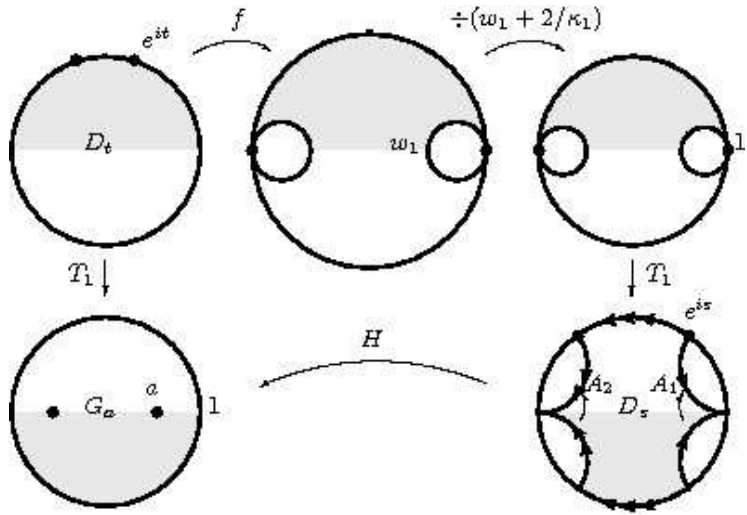}{4.6,-6.5}{6.5}{scale=.7}
 \caption{\small \label{fig:doublepuncture} Universal covering mapping
of doubly-punctured disk.}
\end{figure}

The mapping scheme is depicted in Figure \ref{fig:doublepuncture}. Write  
\[ T_1(z) = \frac{iz+1}{z+i}.
\]
Then $T_1$ sends the upper half of $\D$ to its lower half, and
when 
\begin{equation}\label{eq:tfroma}
   t = \sin^{-1}\left(\frac{1-a^2}{1+a^2}\right),
\end{equation}
$0<t<\pi/2$, we have $T_1(e^{it})=a$.  It can thus be seen that $T_1$  
sends the upper half plane, punctured at the two points $e^{it}$ and
$-e^{-it}$, onto $G_a$.

Now we consider the mapping problem of the disk $\D_t$ to an \scq,
taking $\alpha=0$ and $t$ as just defined.  Let $\lambda_{\rm
  max}$ be the parameter for which $f(e^{it})$ meets $f(e^{-it})$ as
in section \ref{subsec:intervuniv}, where $f=f_{0,t,\lambda_{\rm max}}$.
Recall that this common vertex may be expressed as
$f(e^{it})=w_1-2/\kappa_1$ (note that $\kappa_1=-|\kappa_1|$ is
negative)  so we may divide by this value to obtain a mapping sending
$e^{it}$ to $1$ and sending $1$ to $1/(w_1-2/\kappa_1)$.
 
When we apply $T_1$ to this rescaled \scq, the part in the upper
half-plane is sent to a region in the lower half-plane which may be
described as follows.  It is the lower half of a domain
$D_s\subseteq\D$ which is symmetric in the real axis and bounded by
two arcs of $\partial\D$ ending at $\pm e^{\pm is}$ and four arcs
orthogonal to $\partial\D$. The latter arcs join $1$ to $e^{\pm is}$
and join $-1$ to $-e^{\pm is}$, with tangencies at $1$ and $-1$. By
the above construction, we can deduce that
\begin{equation} \label{eq:sfromwk}
    e^{-is} = T_1(\frac{w_1}{w_1-2/\kappa_1}) .
\end{equation}

  Finally, we define the conformal mapping $H\colon D_s\to G_a$ as the
composition
\[ H(\zeta) = T_1\left(f^{-1}\left(
       (w_1-\frac{2}{\kappa_1})\,T_1^{-1}(\zeta) \right) \right ).
\]

From the above discussion, $H$ sends the upper and lower arcs of
$\partial\D$ bounding the circular hexagon $D_s$ to the upper and
lower semicircle of $\partial\D$ respectively, and sends the four
remaining arcs into the segments $[-1,-a]$ and $[a,1]$.  
 Thus $H$ may be extended by  reflection to all of $\D$, and $D_s$
together with one of its reflected   images forms
a fundamental domain for a Fuchsian group of the second kind
\cite{Mas} acting on $\D$, generated by two M\"obius transformations
$A_1,A_2$ which preserve $\D$ and satisfy  
 \begin{figure}[!t]
\pic{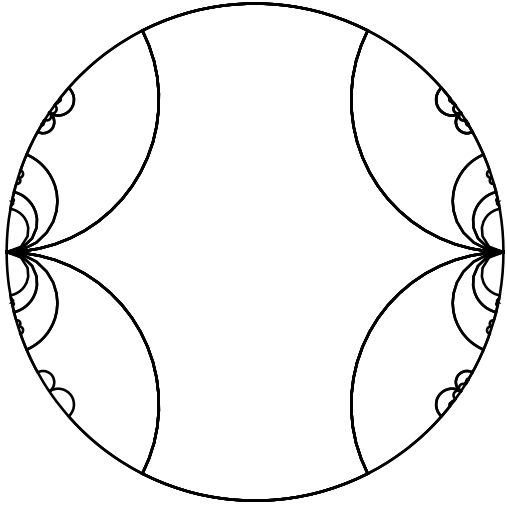}{9.0,0}{6}{scale=.65}
\caption{\small Universal covering  
of doubly-punctured disk.\label{fig:fuchsian}}
\end{figure}
\begin{eqnarray*}
     \mbox{trace}\, A_1 &=& \pm2, \\
     A_1(1) &=&  1, \\
     A_1(e^{-is}) &=& e^{is}, \\
     A_2(-\cg{\zeta}) &=& -\cg{ A_1(\zeta) }.
\end{eqnarray*}

The quotient of $\D$ by this Fuchsian group is conformally equivalent
to $G_a$, and the conformal mapping $H$ (extended to all of $\D$ by
reflection) is the covering map.
 
\begin{table}[th!]
\begin{center}
$\begin{array}{|c|c|c|c|c|}  \hline
\ \ \ a\ \ \ &\ \ \ t\ \ \ & \ \ s/\pi \ \  & \ \  H'(0)\ \ & \mbox{tr}A_1\circ A_2^{-1} \\ \hline\hline
0.1 & \ \;1.37146\ \; & \ 0.32213 \ & 0.00865  & \ \;11.027\\\hline
0.2 & \ \; 1.17601\ \; & 0.26104 & 0.58362 & \ \;19.171\\\hline
0.3 & 0.98788 & 0.21213 & 2.04028 & \ \ 31.388\\\hline
0.4 & 0.80978 & 0.17018 & 3.73469 & \ \ 51.328\\\hline
0.5 & 0.64350 & 0.13330 & 5.31735 & \ \ 86.582\\\hline
0.6 & 0.48996 & 0.10052 & 6.61668 & \ 155.776\\\hline
0.7 & 0.34934 & 0.07123 & 7.57377 & \ 314.856\\\hline
0.8 & 0.22131 & 0.04497 & 8.19717 & \ 797.411\\\hline
0.9 & 0.10517 & 0.02132 & 8.53008 &\!\!\!\! 3560.92\\\hline
\end{array}
$ \vspace{-3ex}
\end{center}
\caption{\small Numerical solution to universal cover of double
  punctured disk.  (About 2.5 seconds of CPU time to compute the
  entire table).
\label{tab:doublepunc} }
\end{table}

The questions posed at the beginning of this section may now be
answered as follows.  Given $a$, we calculate $t$ by
(\ref{eq:tfroma}), and find $\lambda_{\rm max}$ as in section
\ref{subsec:intervuniv}. Based on this parameter we determine
$w_1$,$\kappa_1$ as described in \ref{subsec:p2kappa1}.  This permits
us to obtain $s$ via (\ref{eq:sfromwk}), and then it is
straightforward to find $A_1$ and $A_2$.  Note that $H(0)=0$, because
\[    T_1^{-1}(0) = i,\quad 
     f^{-1}(i(w_1-2/\kappa_1)) = i,\quad
     T_1(i) = 0,
\]
and since
\[   (T_1)'(0)=2,\quad 
    (f^{-1})'(i(w_1-2/\kappa_1))=\frac{1}{f'(i)},\quad 
     (T_1^{-1})'(i)=\frac{1}{2},
\]
the Chain Rule gives us  
\begin{equation}   H'(0) = \frac{w_1-2/\kappa_1}{f'(i)}.
\end{equation}
In the notation of (\ref{eq:fintegral}), $f'(i)=y_1(i)^{-2}$ which we
have seen also can be expressed in terms of the power series produced
by the SPPS method. Thus we may calculate $H'(0)$.
 
 \begin{figure}[!b]
\pic{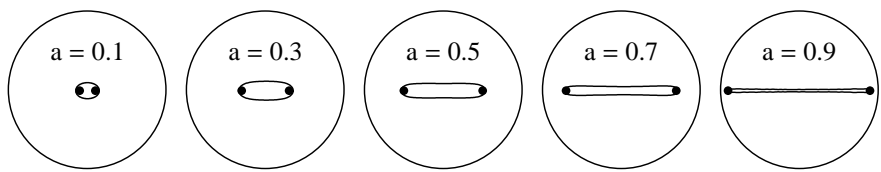}{6.7,0}{2.5}{scale=1.1}
\caption{\label{fig:geodesic} Simple hyperbolic geodesic surrounding the
two punctures of a disk.}
\end{figure}
 
 There is a unique simple geodesic in the hyperbolic metric of $G_a$ which
loops once around the points $\pm a$. This geodesic is covered by the
deck transformation $A_1\circ A_2^{-1}$ of $\D$, which fixes the points  
\begin{equation}
      \frac{-1+e^ {is} \pm \sqrt{2+2e^{2is}}}{ 1+e^{is} }   
\end{equation}
on $\partial\D$. The trace of a hyperbolic M\"obius transformation
serves as a measure of the length of the invariant geodesic which it
covers.  These data can be calculated for any $0<a<1$; sample values
are given in Table \ref{tab:doublepunc}.  It is not difficult to see
that as $a\to1$, $A_1\circ A_2^{-1}$ and the length of the geodesic
tend to $\infty$.  However, $\lim_{a\to1}H'(0)$ apparently is finite.
It would be interesting to see if the values we have calculated could
be obtained as limiting cases of the corresponding values for $\D$
minus two disks (pair of pants) as the radii of these disks tend to
zero.  It would also be interesting to determine the maximum 
of the crossing point of these geodesics with the imaginary axis
(Figure \ref{fig:geodesic}).

Conformal mapping of \scq s can be useful for investigating
many types of  finitely generated Fuchsian groups. In general, 
such groups may define branched coverings rather than true
universal coverings; a ramification of order $n$ corresponds
to a cycle of vertices \cite{Mas} in the fundamental domain whose total
angle is $2\pi/n$. Thus it is of interest to map circle domains with
diverse angles at the vertices.

\noindent Philip R. Brown\\ 
Department of General Academics\\
Texas A\&M University at Galveston\\ 
PO Box 1675, Galveston, Texas 77553 -1675 \\ 
\texttt{ brownp@tamug.edu }

\medskip
\noindent R. Michael Porter \\
Departamento de Matem\'aticas, CINVESTAV--I.P.N.\\
Apdo.\ Postal 1-798, Arteaga 5 \\
Santiago de Queretaro, Qro., 76001  MEXICO \\
\texttt{mike@math.cinvestav.mx}


 \end{document}